\documentclass[11pt,reqno]{amsart}
\usepackage{amssymb,graphicx}

\usepackage{amsfonts,amscd,amsthm,amsmath}

\usepackage{cite}

\numberwithin{equation}{section}

\begin{document}

\begin{center}
\textbf{\large{A dilogarithmic integral arising in quantum field theory}}

\medskip
\textbf{Djurdje Cvijovi\'{c}}
\medskip

{\it Atomic Physics Laboratory, Vin\v{c}a Institute of
Nuclear Sciences \\
P.O. Box $522,$ $11001$ Belgrade$,$ Republic of Serbia}\\

\textbf{E-Mail: djurdje@vinca.rs}\\

\medskip
\textbf{Abstract}

\end{center}

\begin{quotation}

Recently, an interesting dilogarithmic integral arising in quantum field theory has been closed-form evaluated in terms of the Clausen function $\text{Cl}_2(\theta)$  by Coffey [{\it J. Math. Phys.} {\bf 49}  (2008),  093508]. It represents the volume of an ideal tetrahedron in hyperbolic space and is involved in two intriguing equivalent conjectures of Borwein and Broadhurst. It is shown here, by simple and direct arguments, that this integral can be expressed by the triplet of $\,$  the $\,$  Clausen function values which are involved in one of the two above-mentioned conjectures.

\end{quotation}



\vskip .2cm

\section{\bf Introduction and statement of the results}

In a recent paper Coffey \cite[Theorem 1]{Coffey2} has evaluated the integral

\begin{equation} I_{7}= \frac{24}{7 \sqrt7}\int_{\pi/3}^{\pi/2}\ln\left|\frac{\tan(\theta) + \sqrt7}{\tan(\theta)-\sqrt7}\right|d\theta
\end{equation}
\noindent  in a closed form in terms of the Clausen function $\text{Cl}_2(\theta)$ (see Sec. 2)
\begin{equation}
 I_{7}= \frac{12}{7 \sqrt7}\left[\text{Cl}_2(2 \omega_{+})+ 2 \,\text{Cl}_2(\theta_{+})-\text{Cl}_2(2 \omega_{+}+ 2 \theta_{+})\right]\tag{1.2a}
\end{equation}

\noindent where

\begin{equation}\textstyle \omega_{+}=- \textup{tan}^{-1}\left( \frac{2\sqrt3+\sqrt7}{5}\right)
=\textup{tan}^{-1}\left(\sqrt7\right)-\frac{2 \pi}{3}
\quad\textup{and}\quad\theta_{+}=\textup{tan}^{-1}\left(\frac{\sqrt7}{3}\right).\tag{1.2b}
\end{equation}

This and other numerous related integrals arose out of some studies in quantum field theory, in analysis of hyperbolic manifolds whose complementary volumes result from evaluations of Feynman diagrams \cite{Borwein, Broadhurst,Broadhurst1, Broadhurst2,Coffey2,Coffey1,Lunev}. The integral $I_{7}$ is the simplest of $998$ {\em empirically} determined cases where the volume of a manifold is rational multiple of values of various Dirichlet  L series (see Sec. 2). $I_{7}$ was numerically  computed with high accuracy (using highly parallel tanh-sinh quadrature) and represents the volume of an ideal tetrahedron in hyperbolic space.

Two intriguing equivalent conjectures concerning $I_{7}$ have arisen. First, it is  conjectured by Borwein and Broadhurst \cite{Broadhurst} that

\begin{equation}
I_{7} \mathop  = \limits^? \text{L}_{-7}(2)\tag{1.3}
\end{equation}

\noindent where (see Sec. 2)

\begin{align}
\text{L}_{-7}(2)& = \sum_{m = \,0}^{\infty} \textstyle \left[\frac{1}{(7m+1)^2}+\frac{1}{(7m+2)^2}-\frac{1}{(7m+3)^2}+\frac{1}{(7m+4)^2}-\frac{1}{(7m+5)^2}-
\frac{1}{(7m+6)^2}\right]\tag{1.4}
\\
\textstyle
&= \textstyle \frac{1}{49}\left[ \zeta\left(2,\frac{1}{7}\right)+\zeta\left(2,\frac{2}{7}\right)-\zeta\left(2,\frac{3}{7}\right)+
\zeta\left(2,\frac{4}{7}\right)-\zeta\left(2,\frac{5}{7}\right)-\zeta\left(2,\frac{6}{7}\right)\right]\tag{1.4'}
\\
&=\textstyle \frac{2}{\sqrt7}\left[\text{Cl}_2\left(\frac{2 \pi}{7}\right)+ \text{Cl}_2\left(\frac{4\pi}{7}\right)-\text{Cl}_2\left(\frac{6\pi}{7}\right)\right]\tag{1.4"}.
\end{align}

It should be noted that equations (1.4) and (1.4') have been previously given by Coffey \cite{Coffey2}. The ? here indicates that numerical verification  (to $20,000$ digit precision) has been performed but that no formal proof of this "identity" is known \cite{Bailey,Broadhurst}. Second, the $1,800$ digit agreement was found (by integer relation detection algorithm)  for the following unexpected relation between six values of $\text{Cl}_2(\theta)$ ({\em c.f.} Eq. (1.4") above) \cite{Broadhurst}
\begin{equation}\text{L}_{-7}(2)\mathop  = \limits^? \frac{4}{7 \sqrt7} \left[3 \text{Cl}_2(2\phi_{7})-3 \text{Cl}_2(4\phi_{7})+\text{Cl}_2(6\phi_{7})\right],\qquad\phi_{7}=\textup{tan}^{-1}\left(\sqrt7\right).\tag{1.5}
\end{equation}

In this sequel to the work of Coffey \cite{Coffey2, Coffey1} the integral $I_{7}$ is  evaluated in a closed form

\begin{equation}I_{7}= \frac{4}{7 \sqrt7} \left[3 \text{Cl}_2(2\phi_{7})-3 \text{Cl}_2(4\phi_{7})+\text{Cl}_2(6\phi_{7})\right]\tag{1.6}
\end{equation}
\noindent by simple and direct arguments and our result, unlike (1.2a), contains one of the two triplets  of the Clausen function values which are involved in the  conjecture (1.5). We have been unable to deduce (1.6) from Coffey's result (1.2a) and vice versa by the usual integration techniques.

\section{\bf Clausen function and Dirichlet L series}

The Clausen function (of order 2) $\text{Cl}_2(\theta)$, sometimes also called the Clausen integral, is a real function for all $\theta\in\mathbb{R}$ and is given by (\cite[p. 111, Eq. (45)]{Sr-Ch}

\begin{align}\text{Cl}_2(\theta)= \sum_{m=\,1}^{\infty} \frac{\sin\left(m\theta\right)}{m^2}=-\int_{0}^{\theta}\ln\left|2 \sin \left(\frac{t}{2}\right)\right|dt \qquad(\theta\in\mathbb{R}).
\end{align}

\noindent The standard texts on the theory of $\text{Cl}_2(\theta)$ (and various related functions) are Lewin's books \cite{Lewin, Lewin2} and several new results can be found in de Doelder \cite{Doelder}, Grosjean \cite{Grosjean} and Coffey \cite{Coffey2,Coffey1}. Some  elementary properties and special values of $\text{Cl}_2(\theta)$ are:

\begin{align*}
&\textup{a)}\quad\text{Cl}_2(-\theta) = - \text{Cl}_2(\theta),
\\
&\textup{b)}\quad\text{Cl}_2(\theta+2 m \pi)= \text{Cl}_2(\theta), \text{for}\; m\in\mathbb{Z},\hskip65mm
\\
&\textup{c)}\quad\text{Cl}_2(\pi+\theta)= - \text{Cl}_2(\pi-\theta)
\\
&\textup{d)}\quad\text{Cl}_2(m \pi) = 0 , \text{for}\; m\in\mathbb{Z}.
\end{align*}

Define Dirichlet L series (modulo $d$) $\text{L}_d(s)$ in the following way

\begin{equation}
\text{L}_{d}(s)=\sum_{n=\,1}^{\infty}\left(\frac{d}{n}\right)\,\frac{1}{n^s}
\end{equation}

\noindent where $\left(\frac{d}{n}\right)$ is the Kronecker symbol (for more details, see, for instance, \cite{Ayoub}).

Symbol $\left(\frac{d}{n}\right)$ only assumes values $1, -1$ and $0$, it is a periodic function with a period of $|d|$  for an admissible $d$ so that $\text{L}_{d}(s)$ can be expressed as the following finite sum

\begin{equation*}
\text{L}_{d}(s)= \frac{1}{|d|^s} \sum_{\ell =\,1}^{|d|-1}\left(\frac{d}{\ell}\right)\,\zeta\left(s,\frac{\ell}{|d|}\right),
\end{equation*}

\noindent where $\zeta(s,a)$ denotes the Hurwitz (or generalized) zeta function

\begin{equation}
\zeta \left(s, a\right)=\sum_{m = 0}^{\infty}\frac{1}{(m +a)^s}\quad (\mathfrak{R}(s)>1;\, a \not\in \left\{0, -1, -2, \ldots\right\}.
\end{equation}

\noindent If $d<0$, then $\left(\frac{d}{n}\right)$ has the Fourier series expansion

\begin{equation*}
\left(\frac{d}{n}\right)=\frac{1}{\sqrt{|d|}}\sum_{\ell=1}^{|d|-1}\left(\frac{d}{\ell}\right)\sin\left(\frac{2 \,\ell n \pi}{|d|} \right).
\end{equation*}

The sequence $\left(\frac{-7}{n}\right)$, $n\in\mathbb{N}$, is periodic in $n$ with a period of $7$, its first seven values are: $1,1,-1,1,-1,-1,0$ and its expansion is

\begin{equation}\left(\frac{-7}{n}\right)= \frac{2}{\sqrt7} \left[\sin\left(\frac{2\pi}{7}\right)+\sin\left(\frac{4\pi}{7}\right)-\sin\left(\frac{6\pi}{7}\right)\right].
\end{equation}

\noindent It is clear then that

\begin{align*}
&\text{L}_{-7}(s) = \textstyle \frac{1}{7^s}\left[ \zeta\left(s,\frac{1}{7}\right)+\zeta\left(s,\frac{2}{7}\right)-\zeta\left(s,\frac{3}{7}\right)+
\zeta\left(s,\frac{4}{7}\right)-\zeta\left(s,\frac{5}{7}\right)-\zeta\left(s,\frac{6}{7}\right)\right]
\\
&=
\sum_{m=\,0}^{\infty} \textstyle \left[\frac{1}{(7m+1)^s}+\frac{1}{(7m+2)^s}-\frac{1}{(7m+3)^s}+\frac{1}{(7m+4)^s}-\frac{1}{(7m+5)^s}-
\frac{1}{(7m+6)^s}\right],
\end{align*}

\noindent so that for $s=2$ we have the values of $\text{L}_{-7}(2)$ given by (1.4) and (1.4'). In addition, the value of $\text{L}_{-7}(2)$ in terms of the Clausen function (1.4") follows from (2.1), (2.2) and (2.4).

\section{\bf Proof of equation (1.6)}

In order to evaluate the integral given by equation (1.6) we shall need the following two results. The Lemma 2 is familiar multiplication formula for the Clausen function $\text{Cl}_2(\theta)$  (see, for instance, \cite[p. 105, Eq. (4.24)]{Lewin2}) and we include it and  its proof for the sake of self-containedness. In addition, both integrals, (3.1) and (3.2), given by Lemma 1 are essentially the same as the integral in (3.5) below, which, in turn, can be found in well-known books by Lewin \cite {Lewin, Lewin2} (see, for instance, \cite[p. 308, Eq. (37)]{Lewin2}).

\medskip
\noindent{\bf Lemma 1.} In terms of the Clausen function $\text{Cl}_2(\theta)$ defined as in (2.1), we have:

\begin{align}
&\textup{a)}\qquad \int_{\phi}^{x}\ln\left(\frac{\tan(\theta) + \tan(\phi)}{\tan(\theta)-\tan(\phi)}\right)d\theta\nonumber
\\
&\hskip 25mm = - \frac{1}{2} \, \text{Cl}_2(2 x + 2 \phi)+ \frac{1}{2} \,\text{Cl}_2(2 x - 2 \phi)+\frac{1}{2}\,\text{Cl}_2(4 \phi);
\\ \nonumber
\\
&\textup{b)}\qquad \int_{x}^{\phi}\ln\left(\frac{\tan(\phi)+\tan(\theta)}{\tan(\phi)-\tan(\theta)}\right)d\theta\nonumber
\\
& \hskip 25mm =  \frac{1}{2} \,\text{Cl}_2(2 x+2 \phi)- \frac{1}{2} \,\text{Cl}_2(2 x - 2 \phi)- \frac{1}{2}\,\text{Cl}_2(4 \phi).
\end{align}

\noindent{\bf Lemma 2} (Multiplication formula){\bf .} For $m \in \mathbb{N}$, we have

\begin{equation*}
\text{Cl}_2(m\, \theta) = m\, \sum_{\ell=\,0}^{m -1} \text{Cl}_2\left(\theta + \ell \,\frac{2 \pi}{m}\right).
\end{equation*}

In particular

\begin{equation}
\frac{1}{2}\,\text{Cl}_2(2\, \theta) = \text{Cl}_2(\theta)+\text{Cl}_2(\pi+\theta)= \text{Cl}_2(\theta)- \text{Cl}_2(\pi-\theta)
\end{equation}
\noindent and
\begin{equation}
\frac{1}{3}\,\text{Cl}_2(3\, \theta) = \text{Cl}_2(\theta) + \text{Cl}_2\left(\theta+\frac{2 \pi}{3}\right)+ \text{Cl}_2\left(\theta -\frac{2 \pi}{3}\right).
\end{equation}

\medskip
To prove Lemma 1, first observe that $\ln\left[\tan(\theta) + \tan(\phi)\right]$ can be rewritten as follows
\begin{align*}
\ln\left[\tan(\theta) + \tan(\phi)\right]& = \ln\left[\frac{\sin(\theta+\phi)}{\cos(\theta) \,\cos(\phi)}\right]
\\
& \hskip-10mm =\ln\left[\sin\left(\theta+\phi\right)\right] -\ln\left[\sin\left(\frac{\pi}{2}-\theta\right)\right] - \ln\left[\cos(\phi)\right]
\\
&\hskip -10mm =\ln\left[2 \sin\left(\frac{2 \theta +2 \phi}{2}\right)\right]  - \ln\left[2 \sin\left(\frac{\pi-2 \theta}{2}\right)\right] -\ln\left[\cos(\phi)\right]
\end{align*}

\noindent and consequently

\begin{align*}
\int\ln\left[\tan(\theta) + \tan(\phi)\right]d\theta & = - \ln\left[\cos(\phi)\right]\,\theta + \frac{1}{2}\int\ln\left[2 \sin\left(\frac{u}{2}\right)\right]d u
\\
& + \frac{1}{2}\int\ln\left[2 \sin\left(\frac{v}{2}\right)\right]d v + C
\end{align*}

\noindent where

\begin{equation*}
u=2 \theta +2 \phi\qquad\textup{and}\qquad v=\pi-2 \theta,
\end{equation*}

\noindent so that the above integral, in view of the definition of $\text{Cl}_2(\theta)$ in (2.1), at once yields

\begin{align}
\int\ln\left[\tan(\theta)+ \tan(\phi)\right]d\theta = & -\ln\left[\cos(\phi)\right]\,\theta\nonumber
 \\
& - \frac{1}{2}\, \text{Cl}_2(2 \theta+ 2\phi) - \frac{1}{2}\, \text{Cl}_2(\pi - 2\theta)+ C.
\end{align}

\noindent Finally, upon noting that $\text{Cl}_2(0)=0$ (see the definition of $\text{Cl}_2(\theta)$ in (2.1)), the required formulae (3.1) and (3.2) are readily available from (3.5).

It should be noted that Coffey \cite[Sec. 6]{Coffey1} has evaluated around dozen new integrals that are expressible in terms of the function $\text{Cl}_2(x)$, examples of which are

\begin{align*}
&\kappa \,\int_{0}^{x}\ln\left[\sin(\kappa \theta) + \sin(\alpha)\right]d\theta
\\
& = \text{Cl}_2(\alpha)-\text{Cl}_2(\kappa x + \alpha)+ \text{Cl}_2(\alpha - \kappa x + \pi)-\text{Cl}_2(\alpha + \pi)-x\,\kappa \ln 2
\end{align*}

\noindent for $\kappa>0$ and $|x|\leq|\alpha|$, and

\begin{align*}
 - \kappa\, \int_{0}^{x}\ln\left|\cos(\alpha)- \cos(\kappa \theta)\right|d\theta
 = \text{Cl}_2(\kappa\, x- \alpha)+\text{Cl}_2(\kappa\, x+ \alpha)+ x\,\kappa \ln 2,
\end{align*}

\noindent for $\alpha\in\mathbb{R}$, but (3.5) as well as (3.1) and (3.2) are not included in his analysis.

To prove Lemma 2 we apply the well-known trigonometric formula

\begin{equation*}
\sum_{\ell =\,0}^{m -1} \sin\left[n \left(\theta + \ell\,\frac{2 \pi}{m}\right)\right] = \left\{\begin{array}{l}
0, \qquad\qquad\,\,\,\,\,\,\, \textup{if}\, m \nmid n
\\
m \sin(n \theta)\qquad\textup{if}\,m \mid n
\end{array}.
\right.
\end{equation*}

\noindent Hence

\begin{align*}
\sum_{\ell =\,0}^{m-1} \text{Cl}_2\left(\theta + \ell \frac{2 \pi}{m}\right) & = \sum_{n =\,1}^{\infty}\sum_{\ell =\,0}^{m-1}\frac{1}{n^2} \sin\left[n \left(\theta + \ell \frac{2 \pi}{m}\right)\right]
\\
& = \sum_{\scriptsize
 \begin{array}{c} n = \,1
\\
m\mid n \end{array}}^{\infty}\frac{m}{n^2} \sin(n \theta) = \sum_{\ell =\ 1}^{\infty} \frac{m}{(m \ell)^2} \sin(m \ell \theta) = \frac{1}{m}\,\text{Cl}_2(m \theta).
\end{align*}

We conclude this section with the proof of (1.6). Indeed, we may make use of Lemma 1  since the integral $I_{7}$ could be decomposed as

\begin{align}
\frac{7 \sqrt 7}{24}\,I_{7} = \int_{\pi/3}^{\phi_{7}}\ln\left(\frac{\tan(\phi_{7}) + \tan(\theta)}{\tan(\phi_{7})-\tan(\theta)}\right)d\theta
+ \int_{\phi_{7}}^{\pi/2}\ln\left(\frac{\tan(\theta) + \tan(\phi_{7})}{\tan(\theta)-\tan(\phi_{7})}\right)d\theta,
\end{align}

\noindent $\phi_{7}= \textup{tan}^{-1}\left(\sqrt7\right)\approx 1.209429202888189$ being the only singularity inside the interval $\{\pi/3,\pi/2\}$, and in this way we arrive at

\begin{align}
\frac{7 \sqrt 7}{24}\,I_{7} = \frac{1}{2} \left[\text{Cl}_2\left(2 \phi_{7} +\frac{2 \pi}{3}\right)+ \text{Cl}_2\left(2 \phi_{7} - \frac{2 \pi}{3}\right) \right]- \text{Cl}_2\left(\pi + 2 \phi_{7}\right).
\end{align}

\noindent On the other hand, by means of the triplication and duplication formulae, (3.3) and (3.4), we have

\[ \text{Cl}_2\left(2 \phi_{7} + \frac{2 \pi}{3}\right)+ \text{Cl}_2\left(2 \phi_{7} - \frac{2 \pi}{3}\right)= \frac{1}{3}\, \text{Cl}_2\left(6 \phi_{7}\right) -\text{Cl}_2\left(2 \phi_7\right)
\]

\noindent and
\[
\text{Cl}_2\left(\pi + 2 \phi_{7}\right)= \frac{1}{2}\, \text{Cl}_2\left(4 \phi_{7}\right)-\text{Cl}_2\left(2 \phi_{7}\right)
\]

\noindent thus from (3.7) it follows

\begin{align}
\frac{7 \sqrt 7}{24}\,I_{7} = \frac{1}{6} \left[3 \text{Cl}_2\left(2 \phi_{7}\right)-3 \text{Cl}_2\left(4 \phi_{7}\right)+ \text{Cl}_2\left(6 \phi_{7}\right)   \right],
\end{align}

\noindent which is the sought result given by (1.6).

\end{document}